\documentclass[12pt]{article}
\usepackage{graphicx} %% Package for Figure
\usepackage{float} %% Package for Float
\usepackage{amssymb}
\usepackage{biblatex}
\usepackage{amsmath}
\usepackage{indentfirst}
\usepackage{mathtools}
\numberwithin{equation}{section}

\newtheorem{lemma}{Lemma}[section]
\newtheorem{remark}{Remark}[section]
\graphicspath{{figure/}}                                 %% Path of Figures
\usepackage[a4paper]{geometry}                           %% Paper size
\begin{document}
\title{\bf $L^{\vec{p}}-L^{\vec{q}}$ Boundedness of Multiparameter Forelli-Rudin Type Operators on Tube Domains Over The Forward Light Cones
\thanks{The work was Supported by NSFC (11971042, 12071035) , the National Key R\&D Program of China (2021YFA1002600) and Natural Science Foundation of Beijing Municipal (No.1252005).\vskip 1mm
Corresponding author's E-mail address: denggt@bnu.edu.cn }}
\author{\small Xin Xia $^{a}$. Guan Tie Deng$^{b}$\\
 \small a.b\,\,School of Mathematical Sciences, Beijing Normal University, Beijing 100875, China
}
\date{}
\maketitle
\noindent{\bf Abstract:} This study investigates necessary and sufficient conditions for the boundedness of Forelli-Rudin type operators on weighted Lebesgue spaces associated with tubular domains over the forward light cone. We establish a complete characterization of the boundedness for two classes of multiparameter Forelli-Rudin type operators from the mixed-norm Lebesgue space $L^{\vec{p}}$ to $L^{\vec{q}}$, in the range $1 \leq  \vec{p} \leq \vec{q} < \infty$. The findings contribute significantly to the analysis of Bergman projection operators in this setting.   \\
\noindent{\bf{Keywords:}} Forelli-Rudin type operators; Weighted Lebesgue spaces; Boundedness; Forward Light Cone.\\
\smallskip
\noindent {\bf Mathematics Subject Classification: 47B34, 47B38, 47G10.}
\section{Introduction and main results}
As a natural generalization of the classical Lebesgue space \( L^p \), the mixed-norm Lebesgue space \( L^{\vec{p}} \) was introduced by Benedek and Panzone [4] in 1961. Since then, various function spaces with mixed norms have been developed and analyzed. Examples include the mixed-norm weak Lebesgue spaces studied by Chen and Sun [5], and the mixed-norm modulation spaces investigated by Cleanthous and Georgiadis [7, 8]. Owing to their finer structure compared to classical spaces, mixed-norm function spaces have attracted considerable attention in harmonic analysis, partial differential equations, and geometric inequalities [3, 6, 9, 25].

In this article, we study two classes of multiparameter Forelli-Rudin type operators acting on weighted mixed-norm Lebesgue spaces. These operators extend the classical Forelli-Rudin operators introduced in [12]. Our main objective is to establish boundedness conditions for such operators in the setting of tubular domains over the forward light cone. To this end, we first introduce relevant definitions and notations.

We fix a positive integer \(n\) throughout this paper and let \(\mathbb{C}^{n}=\mathbb{C}\times\cdots\times\mathbb{C}\) denote the $n$-dimensional complex Euclidean space. For any two points \(z=(z_{1},\ldots,z_{n})\) and \(w=(w_{1},\ldots,w_{n})\) in \(\mathbb{C}^{n}\), we write
\[
\langle z, w\rangle:=z_{1}\overline{w}_{1}+\cdots+z_{n}\overline{w}_{n}
\]
and \(|z|:=\sqrt{\langle z, z\rangle}\).

For \(z\in\mathbb{C}^{n}\), we also use the notation
\[
z=(z',z_{n}),\quad\text{where }z'=(z_{1},\ldots,z_{n - 1})\in\mathbb{C}^{n - 1}\text{ and }z_{n}\in\mathbb{C}^{1}.
\]

Let \(\mathbb{B}_{n}:=\left\{z\in \mathbb{C}^{n}:\vert z\vert< 1\right\}\) be the open unit ball in \(\mathbb{C}^{n}\) and \(S_{n}:=\left\{z\in \mathbb{C}^{n}:\vert z\vert = 1\right\}\) denote its boundary. For \(a,b,c\in \mathbb{R}\), the Forelli - Rudin type operators are defined as

\[T_{a,b,c}f(z)=(1-\vert z\vert^{2})^{a}\int_{\mathbb{B}_{n}}\frac{(1-\vert u\vert^{2})^{b}f(u)d\nu(u)}{(1-\langle z, u\rangle)^{c}}\]

and

\[S_{a,b,c}f(z)=(1-\vert z\vert^{2})^{a}\int_{\mathbb{B}_{n}}\frac{(1-\vert u\vert^{2})^{b}f(u)d\nu(u)}{\vert1-\langle z, u\rangle\vert^{c}},\]
where $d\nu$ is the volume measure on $\mathbb{B}_{n}$, normalized so that $\nu(\mathbb{B}_{n}) = 1$. Also, for any real parameter $\alpha$ we define $d\nu_{\alpha}(z):=(1 - |z|^{2})^{\alpha}d\nu(z)$. The study of these operators in special cases dates back to Stein \([27]\), who established the boundedness of \(T_{0,0,n+1}\) on \(L^p(\mathbb{B}_n)\) for \(1 < p < \infty\). Subsequently, Forelli and Rudin \([16]\) showed in 1974 that \(T_{0,\sigma+it,n+1+\sigma+it}\) is bounded on \(L^p(\mathbb{B}_n)\) if and only if
\[
(\sigma + 1)p > 1,
\]
where \(1 \leq p < \infty\), \(\sigma > -1\), \(t \in \mathbb{R}\), and \(i = \sqrt{-1}\). Kolaski \([18]\) later revisited this operator from the perspective of Bergman projections. In 1991, Zhu \([29]\) derived the necessary and sufficient condition for the boundedness of the general operator \(T_{a,b,c}\) in the one-dimensional case (\(n = 1\)) with \(c = 1 + a + b\), namely
\[
-pa < \lambda + 1 < p(b + 1).
\]
This result has since been extended to higher dimensions by Kures and Zhu \([20]\), who established the following two theorems.

\textbf{Theorem A} Suppose $1 < p < \infty$. Then the following conditions are equivalent:
\begin{enumerate}
    \item[(i)] The operator $T_{a,b,c}$ is bounded on $L^{p}(\mathbb{B}_n,d\nu_{\alpha})$.
     \item[(ii)] The operator $S_{a,b,c}$ is bounded on $L^{p}(\mathbb{B}_n,d\nu_{\alpha})$.
     \item[(iii)] The parameters satisfy
    \[
    \begin{cases}
    -pa < \alpha + 1 < p(b + 1)\\
    c\leq n + 1+a + b.
    \end{cases}
    \]
\end{enumerate}

\textbf{Theorem B} The following conditions are equivalent:
\begin{enumerate}
    \item[(i)] The operator $T_{a,b,c}$ is bounded on $L^{1}(\mathbb{B}_n,d\nu_{\alpha})$.
    \item[(ii)] The operator $S_{a,b,c}$ is bounded on $L^{1}(\mathbb{B}_n,d\nu_{\alpha})$.
    \item[(iii)] The parameters satisfy
    \[
    \begin{cases}
    -a < \alpha + 1 < b + 1\\
    c=n + 1+a + b
    \end{cases}
    \text{ or }
    \begin{cases}
    -a < \alpha + 1\leq b + 1\\
    c< n + 1+a + b.
    \end{cases}
    \]
\end{enumerate}

These two theorems were originally established in [20] under the additional assumption that \( c \) is neither zero nor a negative integer. Recently, Zhao [35] removed this restriction and further extended the results by characterizing the boundedness of both \( T_{a,b,c} \) and \( S_{a,b,c} \) from \( L^p(\mathbb{B}, dv_{\alpha}) \) to \( L^q(\mathbb{B}, dv_{\beta}) \), where \( 1 \leq p \leq q < \infty \). A key insight in [15, 20] lies in identifying the necessary conditions on the parameter \( c \), derived using inclusion relations among holomorphic function spaces and the radial fractional differential operator. Notably, the special case \( c = n + 1 + a + b \) in Theorem \(\textbf{A}\) is well known, for example, [30, Theorem 2.10].

Before stating our main result, we introduce some definitions and notation. Let \(\Omega\) be an arbitrary subset of \(\mathbb{R}^{n}\). The tubular region \(T_{\Omega}\) is defined as
\[T_{\Omega}=\left\{z=x + iy\in\mathbb{C}^{n}:x\in\mathbb{R}^{n},y\in\Omega\right\}.\]

The definition of the forward light cone $\Lambda_n$ in $\mathbb{R}^n$ is
\[
\Lambda_n = \left\{ (y', y_n) \in \mathbb{R}^n : y_n > |y'| \right\}.
\]

It is obvious that the tube domain $T_{\Lambda_n}$ is a convex set in the complex Euclidean space $\mathbb{C}^n$. For any $y \in \Lambda_n$ and $z \in T_{\Lambda_n}$, we define the functions
\[
\mathfrak{g}(y) = y_n^2 - |y'|^2,
\]
\[
Q(z) := z_1^2 + \cdots + z_{n - 1}^2 - z_n^2.
\]

We can see that $Q(iy) = \mathfrak{g}(y)$. Therefore, the tube domain $T_{\Lambda_n}$ can be expressed as
\[
T_{\Lambda_n} = \left\{ z = x + iy \in \mathbb{C}^n : y_n > 0, \, Q(iy) > 0 \right\}.
\]

Multiparameter operators constitute an important research direction in analysis and operator theory. Over the past several decades, substantial progress has been made in understanding the boundedness of multipliers and singular integral operators in multiparameter function spaces. Notable contributions include the work of Ricci and Stein on multiparameter singular integrals and maximal functions [26], the study by M\"{u}ller, Ricci, and Stein of Marcinkiewicz multipliers and multiparameter structures on Heisenberg groups [23, 24], as well as Fefferman and Pipher's results on multiparameter operators and sharp weighted inequalities [11]. In the present work, we specifically investigate the following two classes of multiparameter operators.
$$T_{\vec{a},\vec{b},\vec{c}} f(z,w)
= \mathfrak{g}(\operatorname{Im} z)^{a_1} \mathfrak{g}(\operatorname{Im} \omega)^{a_2}
\int_{T_{\Lambda_n}} \int_{T_{\Lambda_n}}
\frac{\mathfrak{g}(\operatorname{Im} u)^{b_1} \mathfrak{g}(\operatorname{Im} \eta)^{b_2}}
{Q(z-\overline{u})^{c_1} Q(\omega-\overline{\eta})^{c_2}}
f(u,\eta) \, dV(u) \, dV(\eta)$$
and
$$S_{\vec{a},\vec{b},\vec{c}} f(z,w)
= \mathfrak{g}(\operatorname{Im} z)^{a_1} \mathfrak{g}(\operatorname{Im} \omega)^{a_2}
\int_{T_{\Lambda_n}} \int_{T_{\Lambda_n}}
\frac{\mathfrak{g}(\operatorname{Im} u)^{b_1} \mathfrak{g}(\operatorname{Im} \eta)^{b_2}}
{|Q(z-\overline{u})|^{c_1} |Q(\omega-\overline{\eta})|^{c_2}}
f(u,\eta) \, dV(u) \, dV(\eta)$$
here and thereafter, vectors $\vec{a} := (a_1, a_2), \vec{b} := (b_1, b_2), \vec{c} := (c_1, c_2) \in \mathbb{R}^2$. Obviously, these two
integral operators are natural extensions of the operators $T_{a,b,c}$ and $S_{a,b,c}$. In addition, we would
like to clarify that we limit our focus to the two-dimensional case of the parameters $\vec{a}$, $\vec{b}$ and $\vec{c}$
in this article. However, the higher dimensional case can be approached in a similar manner.

Specifically, we demonstrate that the class of weighted mixed-norm Lebesgue spaces provides a suitable frame-work for analyzing the behavior of the Forelli-Rudin type operators $T_{\vec{a},\vec{b},\vec{c}}$ and $S_{\vec{a},\vec{b},\vec{c}}$. Notice that, for any given $\vec{p} := (p_1, p_2) \in [1, \infty]^2$ and $\vec{\alpha} := (\alpha_1, \alpha_2) \in (-1, \infty)^2$, the weighted mixed-norm
Lebesgue space $L_{\vec{\alpha}}^{\vec{p}} := L^{\vec{p}}\bigl(T_{\Lambda_n} \times T_{\Lambda_n}, dv_{\alpha_1} \, dv_{\alpha_2}\bigr)$ is defined to be the set of all measurable functions $f$ such that $$\|{f}\|_{L_{\vec{\alpha}}^{\vec{p}}}
:= \left\{
\int_{\mathbb{B}_n}
\left[
\int_{\mathbb{B}_n}
\lvert f(z, w) \rvert^{p_1} \, dv_{\alpha_1}(z)
\right]^{\frac{p_2}{p_1}}
\, dv_{\alpha_2}(w)
\right\}^{\frac{1}{p_2}}
< \infty$$
with The measure \( dv_{\alpha}(z) = \mathfrak{g}^{\alpha}(y) dv(z) \). Clearly, when \( p_1 = p_2 = p \) and \( \alpha_1 = \alpha_2 = \alpha \), the space \( L_{\vec{\alpha}}^{\vec{p}} \) then goes back to the weighted Lebesgue space \( L_{\alpha}^p \). The weighted Bergman space \( A_\alpha^p\) on a tubular domain is composed of all holomorphic functions in \( L_\alpha^p \).
It is not difficult to verify that \( A_\alpha^p(T_\Omega) \) is a closed subspace of \( L_\alpha^p(T_\Omega) \).

Finally, we make some conventions on notation. We always denote by \( C \) a positive constant which is independent of the main parameters, but it may vary from line to line. The notation \( f \lesssim g \) means \( f \leq Cg \) and, if \( f \lesssim g \lesssim f \), then we write \( f \sim g \). For any \( q \in [1, \infty] \), we denote by \( q' \) its H\"{o}lder conjugate exponent, namely, \( 1/q + 1/q' = 1 \). Similarly, for any \( \vec{q} := (q_1, q_2) \in [1, \infty]^2 \), we denote by \( \vec{q}' := (q_1', q_2') \) its H\"{o}lder conjugate exponent, namely, \( 1/q_1 + 1/q_1' = 1 \) and \( 1/q_2 + 1/q_2' = 1 \).

With the notion of mixed-norm Lebesgue spaces, our main results are stated as follows. In the rest of this article, we always suppose \( \vec{\alpha} := (\alpha_1, \alpha_2) \in (-1, \infty)^2 \) and \( \vec{\beta} := (\beta_1, \beta_2) \in (-1, \infty)^2 \).

\textbf{Theorem 1} Let \( \vec{p} := (p_1, p_2) \) and \( \vec{q} := (q_1, q_2) \) satisfy \( 1 < p_- \leq p_+ \leq q_- < \infty \) with \( p_+ := \max\{p_1, p_2\} \), \( p_- := \min\{p_1, p_2\} \), and \( q_- := \min\{q_1, q_2\} \). If the operator \( T_{\vec{a},\vec{b},\vec{c}} \) is bounded from \( L_{\vec{\alpha}}^{\vec{p}} \) to \( L_{\vec{\beta}}^{\vec{q}} \), then the parameters satisfy for any \( i \in \{1, 2\} \),
\[
\begin{cases}
-q_i a_i < \beta_i + 1, \ \alpha_i + 1 < p_i (b_i + 1), \\
c_i = n + a_i + b_i + \dfrac{n + \beta_i}{q_i} - \dfrac{n + \alpha_i}{p_i}.
\end{cases}
\]

\textbf{Theorem 2} Let \( 1 < p_- \leq p_+ \leq q_- < \infty \). If the parameters satisfy for any \( i \in \{1, 2\} \), $c_{i}>\frac{3n}{2}$ and
\[
\begin{cases}
-q_i a_i < \beta_i + 1,\alpha_i + 1 < p_i (b_i + 1), \\
c_i = n + \alpha_{i}+ b_i + \dfrac{n + \beta_i}{q_i} - \dfrac{n + \alpha_i}{p_i},
\end{cases}
\]
then the operator \( S_{\vec{a},\vec{b},\vec{c}} \) is bounded from \( L_{\vec{\alpha}}^{\vec{p}} \) to \( L_{\vec{\beta}}^{\vec{q}} \).

\textbf{Theorem 3} Let $\vec{p} := (1, p_2)$ and $\vec{q} := (q_1, q_2)$ satisfy $1 < p_2 \leq q_- < \infty$ with $q_- := \min\{q_1, q_2\}$.
\begin{enumerate}
    \item[(i)] If The operator $T_{\vec{a}, \vec{b}, \vec{c}}$ is bounded from $L_{\vec{\alpha}}^{\vec{p}}$ to $L_{\vec{\beta}}^{\vec{q}}$, then
        $$
    \begin{cases}
    - q_i a_i < \beta_i + 1,\ \alpha_1 < b_1,\ c_1 = a_1 +b_{1}-\alpha_{1} + \frac{\beta_1 + n }{q_{1}}\\
    \alpha_2 + 1 < p_2 (b_2 + 1),\ c_2 = n  + a_2 + b_2 + \dfrac{n  + \beta_2}{q_2} - \dfrac{n  + \alpha_2}{p_2}.
    \end{cases}
    $$
    \item[(ii)] The parameters satisfy that, for any $i \in \{1, 2\}$, $c_{i}>\frac{3n}{2}$,
    $$
    \begin{cases}
    - q_i a_i < \beta_i + 1,\ \alpha_1 < b_1,\ c_1 = a_1 +b_{1}-\alpha_{1} + \frac{\beta_1 + n }{q_{1}}\\
    \alpha_2 + 1 < p_2 (b_2 + 1),\ c_2 = n  + a_2 + b_2 + \dfrac{n  + \beta_2}{q_2} - \dfrac{n  + \alpha_2}{p_2},
    \end{cases}
    $$\end{enumerate}
then the operator $S_{\vec{a}, \vec{b}, \vec{c}}$ is bounded from $L_{\vec{\alpha}}^{\vec{p}}$ to $L_{\vec{\beta}}^{\vec{q}}$.

\textbf{Theorem 4} Let $\vec{p} := (p_1, 1)$ and $\vec{q} := (q_1, q_2)$ satisfy $1 < p_1 \leq q_- < \infty$ with $q_- := \min\{q_1, q_2\}$.

\begin{enumerate}
    \item[(i)] If The operator $T_{\vec{\alpha}, \vec{\beta}, \vec{c}}$ is bounded from $L_{\vec{a}}^{\vec{p}}$ to $L_{\vec{b}}^{\vec{q}}$, then
        $$
    \begin{cases}
    \alpha_1 + 1 < p_1 (b_1 + 1),\ c_1 = n + a_1+ b_1+ \dfrac{n  + \beta_1}{q_1} - \dfrac{n  + \alpha_1}{p_1} \\
    - q_i a_i < \beta_i + 1,\ \alpha_2 < b_2,\ c_2 = a_2+b_{2}-\alpha_{2} + \dfrac{n + \beta_2}{q_2}.
    \end{cases}
    $$
    \item[(ii)] The parameters satisfy that, for any $i \in \{1, 2\}$, $c_{i}>\frac{3n}{2},$
    $$
    \begin{cases}
    \alpha_1 + 1 < p_1 (b_1 + 1),\ c_1 = n + a_1+ b_1+ \dfrac{n  + \beta_1}{q_1} - \dfrac{n  + \alpha_1}{p_1} \\
    - q_i a_i < \beta_i + 1,\ \alpha_2 < b_2,\ c_2 = a_2 + \dfrac{n + \beta_2}{q_2},
    \end{cases}
    $$\end{enumerate}
then the operator $S_{\vec{a}, \vec{b}, \vec{c}}$ is bounded from $L_{\vec{\alpha}}^{\vec{p}}$ to $L_{\vec{\beta}}^{\vec{q}}$.

\textbf{Theorem 5} $\text{Let } 1 = p_1 \leq q_1 < \infty, \, 1 = p_2 \leq q_2 < \infty,$

\begin{enumerate}
    \item[(i)]If the operator $T_{\vec{a}, \vec{b}, \vec{c}}$ is bounded from $L_{\vec{\alpha}}^{\vec{p}}$ to $L_{\vec{\beta}}^{\vec{q}}$, then
     $$
    \begin{cases}
    - q_i a_i < \beta_i + 1,\ \alpha_i < b_i \\
    c_i = a_i +b_{i}-\alpha_{i}+ \dfrac{n + \beta_i}{q_i}.
    \end{cases}
    $$
    \item[(ii)] The parameters satisfy that, for any $i \in \{1, 2\}$, $c_{i}>\frac{3n}{2},$
    $$
    \begin{cases}
    - q_i a_i < \beta_i + 1,\ \alpha_i < b_i \\
    c_i = a_i +b_{i}-\alpha_{i}+ \dfrac{n + \beta_i}{q_i},
    \end{cases}
    $$\end{enumerate}
then the operator $S_{\vec{a}, \vec{b}, \vec{c}}$ is bounded from $L_{\vec{\alpha}}^{\vec{1}}$ to $L_{\vec{\beta}}^{\vec{q}}$.

\section{Preliminaries}
\begin{lemma}
Suppose that \( r, s > \frac{n - 1}{2} \), \( l> -1 \), and \( r + s - l > \frac{3}{2}n - 1 \). Then
\begin{equation}
\int_{T_{\Lambda_n}} \frac{\mathfrak{g}(\operatorname{Im} u)^{l}}{Q(z - \overline{u})^{r} Q(u - \overline{\xi})^{s}} dV(u) \;=\; C_{1,l,r,s} \, Q(z - \overline{\xi})^{n - r - s + l},
\end{equation}

for all \( z \),  $ \xi$ in $T_{\Lambda_n} $, where $C_{1,l,r,s}$  is a constant in terms of n, r, s, l.\end{lemma}

\begin{remark}
$$
    J_{s,l}(z) = \int_{T_{\Lambda_n}} \frac{\mathfrak{g}(\operatorname{Im} u)^l}{|Q(z - \overline{u})|^s} dV(u) =
    \begin{cases}
    C_{2,l,s} \mathfrak{g}(\operatorname{Im} z)^{n - s + l}, & l > -1 \text{ and } l - s < 1 - \frac{3}{2}n \\
    \infty, & \text{otherwise}
    \end{cases},
$$
  where $C_{2,l,s}$  is a constant in terms of n, s, l.

proof: It is clear that $J_{s,l}(z)$ can be derived from  Lemma 2.1.\hfill $\square$
\end{remark}

In order to prove the conclusion of  this article, we introduce the following four lemmas about the Schur's test. All lemmas are derived from~[14].
\begin{lemma}
 Let \( \vec{\mu} := \mu_1 \times \mu_2 \) and \( \vec{\nu} := \nu_1 \times \nu_2 \) be positive measures on the space \( X \times X \), and, for \( i \in \{1,2\} \), \( K_i \) be nonnegative functions on \( X \times X \). Let \( T \) be an integral operator with kernel \( K := K_1 \cdot K_2 \) defined by setting for any \( (x,y) \in X \times X \),
\[
Tf(x,y) := \int_X \int_X K_1(x,s) K_2(y,t) f(s,t) d\mu_1(s) d\mu_2(t).
\]
Suppose \( \vec{p} := (p_1, p_2) \), \( \vec{q} := (q_1, q_2) \in (1, \infty)^2 \) satisfying \( 1 < p_- \leq p_+ \leq q_- < \infty \), where \( p_+ := \max\{p_1, p_2\} \), \( p_- := \min\{p_1, p_2\} \), and \( q_- := \min\{q_1, q_2\} \). Let \( \gamma_i \) and \( \delta_i \) be real numbers such that \( \gamma_i + \delta_i = 1 \) for \( i \in \{1,2\} \). If there exist two positive functions \( h_1 \) and \( h_2 \) defined on \( X \times X \) with two positive constants \( M_1 \) and \( M_2 \) such that for almost all \( (x,y) \in X \times X \),
\begin{equation}
\begin{aligned}
\int_X \left[ \int_X [K_1(x,s)]^{\gamma_1 p_1'} [K_2(y,t)]^{\gamma_2 p_1'} [h_1(s,t)]^{p_1'} d\mu_1(s) \right]^{p_2'/p_1'} d\mu_2(t) \leq M_1 [h_2(x,y)]^{p_2'}
\end{aligned}
\end{equation}
and for almost all \( (s,t) \in X \times X \),
\begin{equation}
\begin{aligned}
\int_X \left[ \int_X [K_1(x,s)]^{\delta_1 q_1} [K_2(y,t)]^{\delta_2 q_1} [h_2(x,y)]^{q_1} d\nu_1(x) \right]^{q_2/q_1} d\nu_2(y) \leq M_2 [h_1(s,t)]^{q_2}, \end{aligned}
\end{equation}
then \( T : L_{\vec{\mu}}^{\vec{p}} \to L_{\vec{\nu}}^{\vec{q}} \) is bounded with \( \| T \|_{L_{\vec{\mu}}^{\vec{p}} \to L_{\vec{\nu}}^{\vec{q}}} \leq M_1^{1/p_2'} M_2^{1/q_2} \).
\end{lemma}

\begin{lemma}
 Let \( \vec{\mu} \), \( \vec{\nu} \), the kernel \( K \), and the operator \( T \) be as in Lemma 2.2. Suppose \( \vec{q} := (q_1, q_2) \in [1, \infty)^2 \). Let \( \gamma_i \) and \( \delta_i \) be two real numbers such that \( \gamma_i + \delta_i = 1 \) for \( i \in \{1, 2\} \). If there exist two positive functions \( h_1 \) and \( h_2 \) defined on \( X \times X \) with two positive constants \( M_1 \) and \( M_2 \) such that for almost all \( (x, y) \in X \times X \),
\begin{equation}
\begin{aligned}
\underset{(s, t) \in X \times X}{\mathrm{ess \, sup}} [K_1(x, s)]^{\gamma_1} [K_2(y, t)]^{\gamma_2} h_1(s, t) \leq M_1 h_2(x, y) \end{aligned}
\end{equation}
and, for almost all \( (s, t) \in X \times X \),
\begin{equation}
\begin{aligned}
\int_X \left[ \int_X [K_1(x, s)]^{\delta_1 q_1} [K_2(y, t)]^{\delta_2 q_1} [h_2(x, y)]^{q_1} d\nu_1(x) \right]^{q_2 / q_1} d\nu_2(y) \leq M_2 [h_1(s, t)]^{q_2}, \end{aligned}
\end{equation}
then \( T : L_{\vec{\mu}}^{\vec{1}} \to L_{\vec{\nu}}^{\vec{q}} \) is bounded with \( \| T \|_{L_{\vec{\mu}}^{\vec{1}} \to L_{\vec{\nu}}^{\vec{q}}} \leq M_1 M_2^{1/q_2} \).
\end{lemma}

\begin{lemma} Let $\vec{\mu}, \vec{v}$, the kernel $K$, and the operator $T$ be as in Lemma 2.2. Suppose $\vec{p} = (p_1, 1)$ with $p_1 \in (1, \infty)$ and $\vec{q} := (q_1, q_2) \in (1, \infty) \times (1, \infty)$ satisfying $1 < p_1 \leq q_- < \infty$, where $q_- := \min\{q_1, q_2\}$. Let $\gamma_i$ and $\delta_i$ be two real numbers such that $\gamma_i + \delta_i = 1$ for $i \in \{1, 2\}$. If there exist two positive functions \( h_1 \) and \( h_2 \) defined on \( X \times X \) with two positive constants \( C_1 \) and \( C_2 \) such that for almost all \( (x, y) \in X \times X \),
\begin{equation}
\begin{aligned}
\underset{t \in X}{\mathrm{ess\ sup}} \int_X (K_1(x, s))^{\gamma_1 p_1'} (K_2(y, t))^{\gamma_2 p_1'} (h_1(s, t))^{p_1'} d\mu_1(s) \leq C_1 (h_2(x, y))^{p_1'}
\end{aligned}
\end{equation}
and, for almost all \( (s, t) \in X \times X \),
\begin{equation}
\begin{aligned}
\int_X \left( \int_X (K_1(x, s))^{\delta_1 q_1} (K_2(y, t))^{\delta_2 q_1} (h_2(x, y))^{q_1} dv_1(x) \right)^{q_2/q_1} dv_2(y) \leq C_2 (h_1(s, t))^{q_2}
\end{aligned}
\end{equation}
then \( T : L_{\vec{\mu}}^{\vec{p}} \to L_{\vec{v}}^{\vec{q}} \) is bounded with \( \| T \|_{L_{\vec{\mu}}^{\vec{p}} \to L_{\vec{v}}^{\vec{q}}} \leq C_1^{1/p_1'} C_2^{1/q_2} \).
\end{lemma}

\begin{lemma} Let \( \vec{\mu}, \vec{v} \), the kernel \( K \), and the operator \( T \) be as in Lemma 2.2. Suppose \( \vec{p} = (1, p_2) \) with \( p_2 \in (1, \infty) \) and \( \vec{q} := (q_1, q_2) \in (1, \infty) \times (1, \infty) \) satisfying \( 1 < p_2 \leq q_- < \infty \), where \( q_- := \min\{q_1, q_2\} \). Let \( \gamma_i \) and \( \delta_i \) be two real numbers such that \( \gamma_i + \delta_i = 1 \) for \( i \in \{1, 2\} \). If there exist two positive functions \( h_1 \) and \( h_2 \) defined on \( X \times X \) with two positive constants \( C_1 \) and \( C_2 \) such that for almost all \( (x, y) \in X \times X \),
\begin{equation}
\begin{aligned}
\int_X \left( \underset{s \in X}{\mathrm{ess\ sup}} (K_1(x, s))^{\gamma_1} (K_2(y, t))^{\gamma_2} h_1(s, t) \right)^{p_2'} d\mu_2(t) \leq C_1 (h_2(x, y))^{p_2'} \end{aligned}
\end{equation}
and, for almost all \( (s, t) \in X \times X \),
\begin{equation}
\begin{aligned}
\int_X \left( \int_X (K_1(x, s))^{\delta_1 q_1} (K_2(y, t))^{\delta_2 q_1} (h_2(x, y))^{q_1} dv_1(x) \right)^{q_2/q_1} dv_2(y) \leq C_2 (h_1(s, t))^{q_2} \end{aligned}
\end{equation}
then \( T : L_{\vec{\mu}}^{\vec{p}} \to L_{\vec{v}}^{\vec{q}} \) is bounded with \( \| T \|_{L_{\vec{\mu}}^{\vec{p}} \to L_{\vec{v}}^{\vec{q}}} \leq C_1^{1/p_2'} C_2^{1/q_2} \).
\end{lemma}

To show the necessity for the boundedness of \( T_{\vec{a}, \vec{b}, \vec{c}} \), we also need the following results.

\begin{lemma}
Suppose $(X, d\mu)$ is a $\sigma$-finite measure space, $1 \leq p < \infty$, and $1/p + 1/q = 1$. Let $G$ be a complex-valued function defined on $X \times X$ and $T$ be the integral operator defined by
\[
Tf(x) = \int_X G(x, y) f(y) d\mu(y).
\]
If the operator $T$ is bounded on $L^p(X, d\mu)$, then its adjoint $T^*$ is the integral operator
\[
T^*f(x) = \int_X \overline{G(y, x)} f(y) d\mu(y)
\]
on $L^q(X, d\mu)$.
\end{lemma}

The following lemma  is a further conclusion on the basis of lemma 2.6.

\begin{lemma}
Let $\vec{p} := (p_1, p_2) \in [1, \infty]^2$ and $\vec{q} := (q_1, q_2) \in [1, \infty]^2$. If the integral operator $T_{\vec{a}, \vec{b}, \vec{c}}$ is bounded from $L_{\vec{\alpha}}^{\vec{p}}$ to $L_{\vec{\beta}}^{\vec{q}}$, then its adjoint operator $T_{\vec{a}, \vec{b}, \vec{c}}^*$ defined by setting
\[
T_{\vec{a}, \vec{b}, \vec{c}}^* g(z, w) := \mathfrak{g}(\operatorname{Im} z)^{b_1 - \alpha_1} \mathfrak{g}(\operatorname{Im} \omega)^{b_2 - \alpha_2} \int_{T_{\Lambda_n}} \int_{T_{\Lambda_n}} \frac{\mathfrak{g}(\operatorname{Im} u)^{a_{1}+\beta_1} \mathfrak{g}(\operatorname{Im} \eta)^{a_{2}+\beta_2}}{Q(z-\overline{u})^{c_1} Q(\omega-\overline{\eta})^{c_2}} g(u, \eta) \, dV(u) \, dV(\eta)
\]
is bounded from $L_{\vec{\beta}}^{\vec{q}'}$ to $L_{\vec{\alpha}}^{\vec{p}'}$.
\end{lemma}

\noindent \textbf{Proof}: Let $g \in L_{\vec{\beta}}^{\vec{q}'}$, we get the definition of operator $T_{\vec{a}, \vec{b}, \vec{c}}^*$ from Lemma 2.6, then
\[
\begin{aligned}
\| T_{\vec{a}, \vec{b}, \vec{c}}^* g \|_{L_{\vec{\alpha}}^{\vec{p}'}} &= \sup_{\| f \|_{L_{\vec{\alpha}}^{\vec{p}}} = 1} \left| \int_{T_{\Lambda_n}} \int_{T_{\Lambda_n}} \overline{f(z, w)} T_{\vec{a}, \vec{b}, \vec{c}}^* g(z, w) \, dV_{\alpha_1}(z) \, dV_{\alpha_2}(w) \right| \\
&= \sup_{\| f \|_{L_{\vec{\alpha}}^{\vec{p}}} = 1} \left| \int_{T_{\Lambda_n}} \int_{T_{\Lambda_n}} \overline{f(z, w)} \mathfrak{g}(\operatorname{Im} z)^{b_1 - \alpha_1} \mathfrak{g}(\operatorname{Im} w)^{b_2 - \alpha_2} \right. \\
&\quad \times \left. \int_{T_{\Lambda_n}} \int_{T_{\Lambda_n}} \frac{\mathfrak{g}(\operatorname{Im} u)^{a_{1}+\beta_1} \mathfrak{g}(\operatorname{Im} \eta)^{a_{2}+\beta_2}}{Q(z-\overline{u})^{c_1} Q(\omega-\overline{\eta})^{c_2}} g(u, \eta) \, dV(u) \, dV(\eta) \, dV_{\alpha_1}(z) \, dV_{\alpha_2}(w) \right|.
\end{aligned}
\]

By the Fubini theorem and the H\"{o}lder inequality of the mixed norms that
\[
\begin{aligned}
\| T_{\vec{a}, \vec{b}, \vec{c}}^* g \|_{L_{\vec{\alpha}}^{\vec{p}'}} &= \sup_{\| f \|_{L_{\vec{\alpha}}^{\vec{p}}} = 1} \left| \int_{T_{\Lambda_n}} \int_{T_{\Lambda_n}} g(u, \eta) \overline{T_{\vec{a}, \vec{b}, \vec{c}} f(u, \eta)} \, dV_{\beta_1}(u) dV_{\beta_2}(\eta) \right| \\
&\leq \sup_{\| f \|_{L_{\vec{\alpha}}^{\vec{p}}} = 1} \left\| T_{\vec{a}, \vec{b}, \vec{c}} f \right\|_{L_{\vec{\beta}}^{\vec{q}}} \| g \|_{L_{\vec{\beta}}^{\vec{q}'}} \lesssim\| g \|_{L_{\vec{\beta}}^{\vec{q}'}}.
\end{aligned}
\]

Thus, the operator  $T_{\vec{a}, \vec{b}, \vec{c}}^*$ is bounded from $L_{\vec{\beta}}^{\vec{q}'}$ to $L_{\vec{\alpha}}^{\vec{p}'}$.
 \hfill $\square$

\section{Proof of Theorem 1}
 We consider the test function \( f_R(z,w) \) defined by
\begin{equation}
f_R(z,w) = \frac{\mathfrak{g}(\operatorname{Im} z)^{l_{1}}\mathfrak{g}(\operatorname{Im} w)^{l_{2}}}{Q(z + \text{iR})^{s_{1}}Q(w + \text{iR})^{s_{2}}} \quad z,w \in T_{\Lambda_n},
\end{equation}
where \( r > 0 \), \( R = (0', r) \in \mathbb{R}^{n - 1} \times \mathbb{R} \), and the real parameters \( s_{i}, l_{i} \) satisfy the conditions
\begin{equation}
\begin{cases}
s_{i} >  \max\left\{ \frac{n}{2}-1, \frac{n-1}{p_{i}} \right\}\\
l_{i} > \max\left\{ -\frac{1 + \alpha_{i}}{p_{i}}, -1 - b_{i} \right\} \\
s_{i} - l_{i} > \max\left\{ \frac{\alpha - 1}{p_{i}}+\frac{3n}{2p}, \frac{3n}{2} - 1 - c_{i} + b_{i} \right\}
\end{cases}.
\end{equation}
We first compute the norm of \( f_R(z,w) \) in $L_{\overrightarrow{\alpha}}^{\overrightarrow{p}}$. Using Lemma 2.1, we deduce that
\begin{equation}
\begin{aligned}
||f_R(z,w)||_{L_{\overrightarrow{\alpha}}^{\overrightarrow{p}}}
&=\left\{
\int_{T_{\Lambda_n}}
\left[
\int_{T_{\Lambda_n}}
\lvert f(z, w) \rvert^{p_1} \, dv_{\alpha_1}(z)
\right]^{\frac{p_2}{p_1}}
\, dv_{\alpha_2}(w)
\right\}^{\frac{1}{p_2}}\\
&=\left\{
\int_{T_{\Lambda_n}}
\left[
\int_{T_{\Lambda_n}}
\lvert \frac{\mathfrak{g}(\operatorname{Im} z)^{l_{1}}\mathfrak{g}(\operatorname{Im} w)^{l_{2}}}{Q(z + \text{iR})^{s_{1}}Q(w + \text{iR})^{s_{2}}} \rvert^{p_1} \, dv_{\alpha_1}(z)
\right]^{\frac{p_2}{p_1}}
\, dv_{\alpha_2}(w)
\right\}^{\frac{1}{p_2}}\\
&=\left( \int_{T_{\Lambda_n}} \frac{\mathfrak{g}(\operatorname{Im} z)^{l_{1}p_{1}+\alpha_{1}}}{|Q(z + {iR})|^{s_{1}p_{1}} } dV(z) \right)^{\frac{1}{p_{1}}}\cdot\left( \int_{T_{\Lambda_n}} \frac{\mathfrak{g}(\operatorname{Im} w)^{l_{2}p_{2}+\alpha_{1}}}{|Q(w + {iR})|^{s_{2}p_{2}} } dV(w) \right)^{\frac{1}{p_{2}}}\\
&=CR^{l_{i}-s_{i}+\frac{n+\alpha_{i}}{p_{i}}}.
\end{aligned}
\end{equation}

where i = 1, 2 and \( C \) is a constant that depends on the parameters \( \alpha_{i}, p_{i},  l_{i}, s_{i} \). The condition (3.2) guarantees that the function \( f_R(z,w) \) belongs to $L_{\vec{\alpha}}^{\vec{p}}$.

We next calculate the norm of \( Tf_R(z,w) \). By Lemma 2.1 and the condition (3.2), we have
$$\begin{aligned}
T_{\vec{a},\vec{b},\vec{c}} f(z,w)
&=
\int_{T_{\Lambda_n}} \int_{T_{\Lambda_n}}
\frac{\mathfrak{g}(\operatorname{Im} z)^{a_1} \mathfrak{g}(\operatorname{Im} \omega)^{a_2}\mathfrak{g}(\operatorname{Im} u)^{b_1+l_{1}} \mathfrak{g}(\operatorname{Im} \eta)^{b_2+l_{2}}}
{Q(z-\overline{u})^{c_1}Q(u+iR)^{s_{1}} Q(\omega-\overline{\eta})^{c_2}Q(\eta+iR)^{s_{2}}}
 dV(u) \, dV(\eta)\\
&=C_{5}\frac{\mathfrak{g}(\operatorname{Im} z)^{a_1} \mathfrak{g}(\operatorname{Im} \omega)^{a_2}}{Q(z+iR)^{c_{1}+s_{1}-n-b_{1}-l_{1}}Q(w+iR)^{c_{2}+s_{2}-n-b_{2}-l_{2}}}.
\end{aligned}
$$

Since the operator \( T \) is bounded on $L_{\vec{\beta}}^{\vec{q}}$, the function \( Tf_R(z,w) \) is in $L_{\vec{\beta}}^{\vec{q}}$. Again by Lemma 2.1 and the condition (3.2), we obtain

\begin{equation}
\begin{cases}
q_{i}a_i + \beta_i > -1 \\
q_{i}(c_{i} - b_{i} - a_{i} - n + s_i - 1_i) - \beta_i >\frac{3n}{2}-1,
\end{cases} \, i = 1, 2.
\end{equation}

Moreover,
$$
\| T f_{\mathbf{R}} \|_{L_{\vec{\beta}}^{\vec{q}}} = C' R^{a_{i} + b_{i} - c_{i} + l_i - s_i + n  + \frac{\beta_i  + n}{q_{i}}},
$$
where i = 1, 2, \( C' \) is a constant depending only on $ a_{i}, b_{i}, c_{i}, l_{i}, s_{i}, \vec{\alpha}, \vec{\beta}$, $p_{i}$ and $q_{i}$. Due to the boundedness of the operator \( T \) from $L_{\vec{\alpha}}^{\vec{p}}$ to $L_{\vec{\beta}}^{\vec{q}}$, we have
$$
 \| T f_{\mathbf{R}} \|_{L_{\vec{\beta}}^{\vec{q}}} \leq C''\| f_{\mathbf{R}} \|_{L_{\vec{\alpha}}^{\vec{p}}}.
$$
That is,
$$
C'R^{a_{i} + b_{i} - c_{i} + l_i - s_i + n  + \frac{\beta_i  + n}{q_{i}}+ } \leq CC'' R^{l_{i}-s_{i}+\frac{n+\alpha_{i}}{p_{i}}}  ,
$$
where \( C, C' \) and \( C'' \) are independent of $R$. For the selection of $R$, it is only required to be positive integers. Therefore, for the above formula to hold, the following condition must be satisfied:
\begin{equation}
c_i = a_i + b_i + n  + \frac{\beta_1 + n }{q_{i}} - \frac{\alpha_1 + n }{p_{i}}, \quad i = 1, 2.
\end{equation}
Combining conditions (3.2) and (3.5), condition (3.4) is equivalent to
\begin{equation}
-a_i q_{i} < \beta_i + 1, \quad i= 1, 2.
\end{equation}

From Lemma 2.7, $T_{\vec{a}, \vec{b}, \vec{c}}$ is bounded from $L_{\vec{\alpha}}^{\vec{p}}$ to $L_{\vec{\beta}}^{\vec{q}}$, then its adjoint operator $T_{\vec{a}, \vec{b}, \vec{c}}^*$ is bounded from $L_{\vec{\beta}}^{\vec{q}'}$ to $L_{\vec{\alpha}}^{\vec{p}'}$.
Using the same discussion as above we derive that
\[
-p_{i}'(b_{i} - \alpha_{i}) < \alpha_{i} + 1 ,
\]
which implies
\begin{equation}
\alpha_{i} + 1 < p_{i}(b_{i} + 1) .
\end{equation}

 \hfill $\square$

\section{Proof of Theorem 2}
For any \( i \in \{1, 2\} \), define
\[
\lambda_i := \frac{n  + \beta_i}{q_i} - \frac{n + \alpha_i}{p_i}, \ c_i := n + a_i +b_i + \lambda_i, \ \text{and} \ \tau_i := c_i -a_i - b_i + \alpha_i.
\]
By the fact \( -(1 + \beta_i)/q_i < 0 \), we know that there exist two negative numbers \( r_1 \) and \( r_2 \) such that, for any \( i \in \{1, 2\} \), \( -\frac{1 + \beta_i}{q_i} < r_i < 0 \). In addition, from \( \alpha_i + 1 < p_i(b_i + 1) \), it follows that, for any \( i \in \{1, 2\} \),
\begin{equation}
\begin{aligned}
b_i - \alpha_i + \frac{\alpha_i + 1}{p_i'} > 0.
\end{aligned}
\end{equation}

Notice that, for any \( i \in \{1, 2\} \), \( \tau_i = \frac{n + \alpha_i}{p_i'} + \frac{n + \beta_i}{q_i} > 0 \), which, combined with (4.1), further implies that
\[
-\frac{\tau_i(1 + \alpha_i)}{p_i'} - \frac{(b_i - \alpha_i)(n + \alpha_i)}{p_i'} < \frac{(b_i - \alpha_i)(n + \beta_i)}{q_i}.
\]
Then there exist \( s_1 \) and \( s_2 \) such that, for any \( i \in \{1, 2\} \),
\[
-\frac{\tau_i(1 + \alpha_i)}{p_i'} - \frac{(b_i - \alpha_i)(n + \alpha_i)}{p_i'} < \tau_i s_i + (b_i - \alpha_i)(s_i - r_i) < \frac{(b_i - \alpha_i)(n + \beta_i)}{q_i},
\]
which is equivalent to
\begin{equation}
\begin{aligned}
-\frac{1 + \alpha_i}{p_i'} - (b_i - \alpha_i)\gamma_i < s_i < (b_i - \alpha_i)\delta_i, \end{aligned}
\end{equation}
where, for any \( i \in \{1, 2\} \), \( \gamma_i := \frac{(n+ \alpha_i)/p_i' + s_i - r_i}{\tau_i} \) and \( \delta_i := \frac{(n + \beta_i)/q_i + r_i - s_i}{\tau_i} \).

Obviously, \( \gamma_i + \delta_i = 1 \). Let \( h_1(u, \eta) := \mathfrak{g}(\operatorname{Im} u)^{s_1} \mathfrak{g}(\operatorname{Im} \eta)^{s_2} \), \( h_2(z, w) := \mathfrak{g}(\operatorname{Im} z)^{r_1} \mathfrak{g}(\operatorname{Im} \omega)^{r_2} \),
\[
K_1(z, u) := \frac{\mathfrak{g}(\operatorname{Im} z)^{a_{1}}\mathfrak{g}(\operatorname{Im} u))^{b_1 - \alpha_1}}{|Q(z-\overline{u})|^{c_1}}, \ \text{and} \ K_2(w, \eta) = \frac{\mathfrak{g}(\operatorname{Im} w)^{a_{2}}\mathfrak{g}(\operatorname{Im} \eta)^{b_2 - \alpha_2}}{|Q(\omega-\overline{\eta})|^{c_2}}.
\]

In order to use lemma 2.2 to show the boundedness of \( S_{\vec{a}, \vec{b}, \vec{c}} \), we consider
\begin{equation}
\begin{aligned}
& \int_{T_{\Lambda_n}} [K_1(z, u)]^{\gamma_1 p_1'} [K_2(w, \eta)]^{\gamma_2 p_1'} [h_1(u, \eta)]^{p_1'} dV_{\alpha_1}(u) \\
& = \frac{\mathfrak{g}(\operatorname{Im} w)^{a_{2}\gamma_2 p_1'}\mathfrak{g}(\operatorname{Im} \eta)^{(b_2 - \alpha_2)\gamma_2 p_1' + s_2 p_1'}}{|Q(\omega-\overline{\eta})|^{c_2 \gamma_2 p_1'}} \\ &\cdot\int_{T_{\Lambda_n}} \frac{\mathfrak{g}(\operatorname{Im} z)^{a_{1}\gamma_1 p_1'}\mathfrak{g}(\operatorname{Im} u)^{(b_1 - \alpha_1)\gamma_1 p_1' + s_1 p_1' + \alpha_1}}{|Q(z-\overline{u})|^{c_1 \gamma_1 p_1'}} dV(u).
\end{aligned}
\end{equation}
From the left inequality of (4.2), we have
\begin{equation}
\begin{aligned}
(b_i - \alpha_i)\gamma_i p_i' + s_i p_i' + \alpha_i > -1.
\end{aligned}
\end{equation}
Moreover, $$(b_i - \alpha_i)\gamma_i p_i' + s_i p_i' + \alpha_i -c_1 \gamma_1 p_1'< 1-\frac{3n}{2},$$
by the fact that \( (c_i - b_i-a_{i} + \alpha_i)\gamma_i = \tau_i \gamma_i = \frac{n  + \alpha_i}{p_i'} + s_i - r_i \), it follows that, for any \( i \in \{1, 2\} \),
\begin{equation}
\begin{aligned}
n+ (a_{i}+b_i - \alpha_i)\gamma_i p_i' + s_i p_i' + \alpha_i - c_i \gamma_i p_i' = r_i p_i'. \end{aligned}
\end{equation}
From this, (4.4), and Remark 2.1, we infer that, for any given \( z \in T_{\Lambda_n} \), we have
\[
\int_{T_{\Lambda_n}} \frac{\mathfrak{g}(\operatorname{Im} z)^{a_{1}\gamma_1 p_1'}\mathfrak{g}(\operatorname{Im} u)^{(b_1 - \alpha_1)\gamma_1 p_1' + s_1 p_1' + \alpha_1}}{|Q(z-\overline{u})|^{c_1 \gamma_1 p_1'}} dV(u) \lesssim \mathfrak{g}(\operatorname{Im} z))^{r_1 p_1'},
\]
which, together with (4.3), (4.4), (4.5), and Lemma 2.2, we can get
\[
\begin{aligned}
& \int_{T_{\Lambda_n}} \left[ \int_{T_{\Lambda_n}} [K_1(z, u)]^{\gamma_1 p_1'} [K_2(w, \eta)]^{\gamma_2 p_1'} [h_1(u, \eta)]^{p_1'} dV_{\alpha_1}(u) \right]^{p_2'/p_1'} dV_{\alpha_2}(\eta) \\
& \quad \lesssim \mathfrak{g}(\operatorname{Im} z)^{r_1 p_2'} \int_{T_{\Lambda_n}} \frac{\mathfrak{g}(\operatorname{Im} w)^{a_{2}\gamma_2 p_1'}\mathfrak{g}(\operatorname{Im} \eta)^{(b_2 - \alpha_2)\gamma_2 p_2' + s_2 p_2' + \alpha_2}}{|Q(\omega-\overline{\eta})|^{c_2 \gamma_2 p_2'}} dV(\eta)\\
& \quad \lesssim \mathfrak{g}(\operatorname{Im} z)^{r_1 p_2'} \mathfrak{g}(\operatorname{Im} w)^{r_2 p_2'} = [h_2(z, w)]^{p_2'}.
\end{aligned}
\]
Thus, condition (2.2) holds true for the operator \( S_{\vec{a}, \vec{b}, \vec{c}} \).

Next, We check the second condition (2.3). Notice that
\[
\begin{aligned}
& \int_{T_{\Lambda_n}} [K_1(z, u)]^{\delta_1 q_1} [K_2(w, \eta)]^{\delta_2 q_1} [h_2(z, w)]^{q_1} dV_{\beta_1}(z) \\
& \quad = \frac{\mathfrak{g}(\operatorname{Im} \eta)^{(b_2 - \alpha_2)\delta_2 q_1} \mathfrak{g}(\operatorname{Im} w)^{a_{2}\delta_2 q_1+\gamma_2 q_1}}{|Q(w-\overline{\eta})|^{c_2 \delta_2 q_1}} \\
& \quad \cdot \int_{T_{\Lambda_n}} \frac{\mathfrak{g}(\operatorname{Im} u)^{(b_1 - \alpha_1)\delta_1 q_1} \mathfrak{g}(\operatorname{Im} z)^{a_{1}\delta_1 q_1+\gamma_1 q_1 + \beta_1}}{|Q(z-\overline{u})|^{c_1 \delta_1 q_1}} dV(z).
\end{aligned}
\]
We obviously have
\[
r_i q_i + \beta_i > -1. \tag{4.6}
\]
In addition, from the choice of \( \delta_i \), it follows that, for any \( i \in \{1, 2\} \), \( (c_i -a_{i}- b_i + \alpha_i)\delta_i = \tau_i \delta_i = \frac{n + \beta_i}{q_i} + r_i - s_i \), which, combined with the right inequality of (4.2), we have
\[
n+ r_i q_i + \beta_i - c_i \delta_i q_i = s_i q_i - (a_{i}+b_i - \alpha_i)\delta_i q_i . \tag{4.7}
\]

From (4.6), (4.7)and Remark 2.1, we get
\[
\int_{T_{\Lambda_n}} \frac{\mathfrak{g}(\operatorname{Im}z))^{a_{1}\delta_1 q_1+\gamma_1 q_1 + \beta_1}}{|Q(z-\overline{u})|^{c_1 \delta_1 q_1}} dV(z) \lesssim \mathfrak{g}(\operatorname{Im} u)^{s_1 q_1 - (a_{1}+b_1 - \alpha_1)\delta_1 q_1}.
\]

Then together with (4.6), (4.7), and Lemma 2.2 again, we have
\[
\begin{aligned}
& \int_{T_{\Lambda_n}} \left[ \int_{T_{\Lambda_n}} [K_1(z, u)]^{\delta_1 q_1} [K_2(w, \eta)]^{\delta_2 q_1} [h_2(z, w)]^{q_1} dV_{\beta_1}(z) \right]^{q_2 / q_1} dV_{\beta_2}(w) \\
& \quad \lesssim \mathfrak{g}(\operatorname{Im} u)^{s_1 q_2} \mathfrak{g}(\operatorname{Im} \eta)^{(a_{2}+b_2 - \alpha_2)\delta_2 q_2} \int_{T_{\Lambda_n}} \frac{\mathfrak{g}(\operatorname{Im} w))^{\gamma_2 q_2 + \beta_2}}{|Q(w-\overline{\eta})|^{c_2 \delta_2 q_2}} dV(\eta) \\
& \quad \lesssim \mathfrak{g}(\operatorname{Im} u)^{s_1 q_2} \mathfrak{g}(\operatorname{Im} \eta)^{s_2 q_2} = [h_1(u, \eta)]^{q_2}.
\end{aligned}
\]

Thus, the operator \( S_{\vec{a}, \vec{b}, \vec{c}} \) satisfies all conditions of Lemma 2.2, then \( S_{\vec{a}, \vec{b}, \vec{c}} \) is bounded from \( L_{\vec{\alpha}}^{\vec{p}} \) to \( L_{\vec{\beta}}^{\vec{q}} \).

 \hfill $\square$

\section{Proof of Theorem 3}
(i) We consider the test function \( f_R(z,w) \) defined by
\begin{equation}
f_R(z,w) = \frac{\mathfrak{g}(\operatorname{Im} z)^{l_{1}}}{Q(z + \text{iR})^{s_{1}}Q(w + \text{iR})^{s_{2}}} \quad z,w \in T_{\Lambda_n},
\end{equation}
where \( r > 0 \), \( R = (0', r) \in \mathbb{R}^{n - 1} \times \mathbb{R} \), and the real parameters \( s_{i}, l_{i} \) satisfy the conditions
\begin{equation}
\begin{cases}
s_{i} >  \max\left\{ \frac{n}{2}-1, \frac{n-1}{p_{i}} \right\}\\
l_{i} > \max\left\{ -\frac{1 + \alpha_{i}}{p_{i}}, -1 - b_{i} \right\} \\
s_{i} - l_{i} > \max\left\{ \frac{\alpha - 1}{p_{i}}+\frac{3n}{2p}, \frac{3n}{2} - 1 - c_{i} + b_{i} \right\}
\end{cases}.
\end{equation}
We first compute the norm of \( f_R(z,w) \) in $L_{\overrightarrow{\alpha}}^{\overrightarrow{p}}$. Using Lemma 2.1, we deduce that
\begin{equation}
\begin{aligned}
||f_R(z,w)||_{\overrightarrow{\alpha}}^{\overrightarrow{p}}
&=\left\{
\int_{T_{\Lambda_n}}
\left[
\int_{T_{\Lambda_n}}
\lvert f(z, w) \rvert \, dV_{\alpha_1}(z)
\right]^{p_2}
\, dV_{\alpha_2}(w)
\right\}^{\frac{1}{p_2}}\\
&=\left\{
\int_{T_{\Lambda_n}}
\left[
\int_{T_{\Lambda_n}}
\lvert \frac{\mathfrak{g}(\operatorname{Im} z)^{l_{1}}}{Q(z + \text{iR})^{s_{1}}Q(w + \text{iR})^{s_{2}}} \rvert \, dv_{\alpha_1}(z)
\right]^{p_2}
\, dv_{\alpha_2}(w)
\right\}^{\frac{1}{p_2}}\\
&=\left( \int_{T_{\Lambda_n}} \frac{\mathfrak{g}(\operatorname{Im} z)^{l_{1}+\alpha_{1}}}{|Q(z + {iR})|^{s_{1}} } dV(z) \right)\cdot\left( \int_{T_{\Lambda_n}} \frac{\mathfrak{g}(\operatorname{Im} w)^{\alpha_{2}}}{|Q(w + {iR})|^{s_{2}p_{2}} } dV(w) \right)^{\frac{1}{p_{2}}}\\
&=CR^{n+l_{1}-s_{1}+\alpha_{1}}\cdot R^{-s_{2}+\frac{n+\alpha_{2}}{p_{2}}}.
\end{aligned}
\end{equation}

\( C \) is a constant that depends on the parameters \( \alpha_{i}, p_{i},  l_{i}, s_{i} \). The condition (5.2) guarantees that the function \( f_R(z,w) \) belongs to $L_{\overrightarrow{\alpha}}^{\overrightarrow{p}}$.

We next calculate the norm of \( Tf_R(z,w) \). By Lemma 2.1 and the condition (5.2), we have
$$\begin{aligned}
T_{\vec{a},\vec{b},\vec{c}} f(z,w)
&=
\int_{T_{\Lambda_n}} \int_{T_{\Lambda_n}}
\frac{\mathfrak{g}(\operatorname{Im} z)^{a_1} \mathfrak{g}(\operatorname{Im} \omega)^{a_2}\mathfrak{g}(\operatorname{Im} u)^{b_1+l_{1}}\mathfrak{g}(\operatorname{Im} \eta)^{b_2}}
{Q(z-\overline{u})^{c_1}Q(u+iR)^{s_{1}} Q(\omega-\overline{\eta})^{c_2}Q(\eta+iR)^{s_{2}}}
 dV(u) \, dV(\eta)\\
&=C_{1}\frac{\mathfrak{g}(\operatorname{Im} z)^{a_1} \mathfrak{g}(\operatorname{Im} \omega)^{a_2}}{Q(z+iR)^{c_{1}+s_{1}-n-b_{1}-l_{1}}Q(w+iR)^{c_{2}+s_{2}-n-b_{2}}}.
\end{aligned}
$$

Since the operator \( T \) is bounded on $L_{\overrightarrow{\beta}}^{\overrightarrow{q}}$, the function \( Tf_R(z,w) \) is in $L_{\overrightarrow{\beta}}^{\overrightarrow{q}}$. Again by Lemma 2.1 and the condition (5.2), we obtain

\begin{equation}
\begin{cases}
q_{i}a_i + \beta_i > -1 \\
q_{1}(c_{1} - b_{1} - a_{1} - n + s_1 - l_1) - \beta_1 >\frac{3n}{2}-1\\
q_{2}(c_{2} - b_{2} - a_{2} - n + s_2) - \beta_2 >\frac{3n}{2}-1,
\end{cases} \, i = 1, 2.
\end{equation}

Moreover,
$$
\| T f_{\mathbf{R}} \|_{L_{\overrightarrow{\beta}}^{\overrightarrow{q}}} = C' R^{a_{1} + b_{1} - c_{1} + l_1 - s_1 + n  + \frac{\beta_1  + n}{q_{1}}}\cdot R^{a_{2} + b_{2} - c_{2} - s_2 + n  + \frac{\beta_2  + n}{q_{2}}},
$$
\( C' \) is a constant depending only on $ a_{i}, b_{i}, c_{i}, l_{1}, s_{i}$, $\overrightarrow{\alpha}, \overrightarrow{\beta}$, $p_{2}$ and $q_{i}$, where $i$ = 1, 2. Due to the boundedness of the operator \( T \) from $L_{\overrightarrow{\alpha}}^{\overrightarrow{p}}$ to $L_{\overrightarrow{\beta}}^{\overrightarrow{q}}$, we have
$$
 \| T f_{\mathbf{R}} \|_{L_{\overrightarrow{\beta}}^{\overrightarrow{q}}} \leq C''\| f_{\mathbf{R}} \|_{L_{\overrightarrow{\alpha}}^{\overrightarrow{p}}}.
$$
That is,
$$
\begin{cases}C'R^{a_{1} + b_{1} - c_{1} + l_1 - s_1 + n  + \frac{\beta_1  + n}{q_{1}}+ } \leq CC'' R^{l_{1}-s_{1}+n+\alpha_{1}} \\
C'R^{a_{2} + b_{2} - c_{2} - s_2 + n  + \frac{\beta_2  + n}{q_{2}} } \leq CC'' R^{-s_{i}+\frac{n+\alpha_{2}}{p_{2}}},
\end{cases}
$$
where \( C, C' \) and \( C'' \) are independent of $R$. For the selection of $R$, it is only required to be positive integers. Therefore, for the above formula to hold, the following condition must be satisfied:
\begin{equation}
\begin{cases}
\alpha_{1}<b_{1}\\
c_1 = a_1 +b_{1}-\alpha_{1} + \frac{\beta_1 + n }{q_{1}}, \\
c_2 = a_2 + b_2 + n  + \frac{\beta_2 + n }{q_{2}} - \frac{\alpha_2 + n }{p_{2}}.
\end{cases}
\end{equation}

Combining conditions (5.2) and (5.5), condition (5.4) is equivalent to
\begin{equation}
-a_i q_{i} < \beta_i + 1, \quad i= 1, 2.
\end{equation}

From Lemma 2.7, $T_{\vec{a}, \vec{b}, \vec{c}}$ is bounded from $L_{\vec{\alpha}}^{\vec{p}}$ to $L_{\vec{\beta}}^{\vec{q}}$, then its adjoint operator $T_{\vec{a}, \vec{b}, \vec{c}}^*$ is bounded from $L_{\vec{\beta}}^{\vec{q}'}$ to $L_{\vec{\alpha}}^{\vec{p}'}$.
Using the same discussion as above we derive that
\[
-p_{2}'(b_{2} - \alpha_{2}) < \alpha_{2} + 1 ,
\]
which implies
\begin{equation}
\alpha_{2} + 1 < p(b_{2} + 1) .
\end{equation}

(ii) This case is a symmetric case of Theorem 4(ii), hence we omit the proof.

 \hfill $\square$

\section{Proof of Theorem 4}

(i) This case is a symmetric case of Theorem 3(i). Therefore, its proof is similar to that of Theorem 1 and we omit the proof.

(ii) Proof. Suppose
\[
\lambda_1 := \frac{n  + \beta_1}{q_1} - \frac{n + \alpha_1}{p_1}, \ c_1 := n + a_1 +b_1 + \lambda_1, \ \text{and} \ \tau_1 := c_1 -a_1 - b_1 + \alpha_1.
\]
By the fact \( -(1 + \beta_i)/q_i < 0 \), we know that there exist two negative numbers \( r_1 \) and \( r_2 \) such that, for any \( i \in \{1, 2\} \), \( -\frac{1 + \beta_i}{q_i} < r_i < 0 \). In addition, for any \( i \in \{1, 2\} \),
we have $\tau_1 = \frac{n+1+\alpha_1}{p_1'} + \frac{n+1+\beta_1}{q_1} > 0$. This, combined with the fact that $a_{1}+b_1 - \alpha_1 + \frac{\alpha_1 + 1}{p_1'} > 0$, further implies that,
\[
-\frac{\tau_1(1 + \alpha_1)}{p_1'} - \frac{(a_{1}+b_1 - \alpha_1)(n+ \alpha_1)}{p_1'} < \frac{(a_{1}+b_1 - \alpha_1)(n+ \beta_1)}{q_1}.
\]
Thus, there exists $s_1$ such that
\[
-\frac{\tau_1(1 + \alpha_1)}{p_1'} - \frac{(a_{1}+b_1 - \alpha_1)(n + \alpha_1)}{p_1'} < \tau_1 s_1 + (a_{1}+b_1 - \alpha_1)(s_1 - r_1) < \frac{(a_{1}+b_1 - \alpha_1)(n + \beta_1)}{q_1}.
\]
which is equivalent to
\begin{equation}
\begin{aligned}
-\frac{1 + \alpha_1}{p_1'} - (a_1 +b_1 - \alpha_1)\gamma_1 < s_1 < (a_1 +b_1 - \alpha_1)\delta_1, \end{aligned}
\end{equation}
where \( \gamma_1 := \frac{(n+ \alpha_1)/p_1' + s_1 - r_1}{\tau_1} \) and \( \delta_1 := \frac{(n + \beta_1)/q_1 + r_1 - s_1}{\tau_1} \). Clearly, \( \gamma_1 + \delta_1 = 1 \).
Let $\gamma_2 :=\frac{s_2-r_{2}}{\tau_2}$  and  \( \delta_2:= \frac{(n + \beta_2)/q_2 + r_2 - s_2}{\tau_2} \).  Obviously, \( \gamma_2 + \delta_2 = 1 \).
We now define \( h_1(u, \eta) := \mathfrak{g}(\operatorname{Im} u)^{s_1}\mathfrak{g}(\operatorname{Im} \eta)^{s_2} \), \( h_2(z, w) := \mathfrak{g}(\operatorname{Im} z)^{r_1} \mathfrak{g}(\operatorname{Im} \omega)^{r_2} \),
\[
K_1(z, u) := \frac{\mathfrak{g}(\operatorname{Im} z)^{a_{1}}\mathfrak{g}(\operatorname{Im} u)^{b_1 - \alpha_1}}{|Q(z-\overline{u})|^{c_1}}, \ \text{and} \ K_2(w, \eta) = \frac{\mathfrak{g}(\operatorname{Im} w)^{a_{2}}\mathfrak{g}(\operatorname{Im} \eta)^{b_2 - \alpha_2}}{|Q(\omega-\overline{\eta})|^{c_2}}.
\]

In order to use lemma 2.4, we consider
\begin{equation}
\begin{aligned}
& \int_{T_{\Lambda_n}} [K_1(z, u)]^{\gamma_1 p_1'} [K_2(w, \eta)]^{\gamma_2 p_1'} [h_1(u, \eta)]^{p_1'} dV_{\alpha_1}(u) \\
& = \frac{\mathfrak{g}(\operatorname{Im} w)^{a_{2}\gamma_2 p_1'}\mathfrak{g}(\operatorname{Im} z)^{a_{1}\gamma_1 p_1'}}{|Q(\omega-\overline{\eta})|^{c_2 \gamma_2 p_1'}}\int_{T_{\Lambda_n}} \frac{\mathfrak{g}(\operatorname{Im} u)^{(b_1 - \alpha_1)\gamma_1 p_1' + s_1 p_1' + \alpha_1}}{|Q(z-\overline{u})|^{c_1 \gamma_1 p_1'}} dV(u).
\end{aligned}
\end{equation}
From the left inequality of (6.1), we have
\begin{equation}
\begin{aligned}
(a_{1}+b_1 - \alpha_1)\gamma_i p_1' + s_1 p_1' + \alpha_1 > -1.
\end{aligned}
\end{equation}
Moreover, by the fact that \( (c_1 - b_1-a_{1} + \alpha_1)\gamma_1 = \tau_1 \gamma_1 = \frac{n  + \alpha_1}{p_1'} + s_1 - r_1 \), it follows that,
\begin{equation}
\begin{aligned}
n+ (a_{1}+b_1- \alpha_1)\gamma_1 p_1' + s_1 p_1' + \alpha_1 - c_1 \gamma_1 p_1' = r_1 p_1' < 0. \end{aligned}
\end{equation}
From this, (6.4), and Lemma 2.1, we infer that, for any given \( z \in T_{\Lambda_n}\),
\[
\int_{T_{\Lambda_n}} \frac{\mathfrak{g}(\operatorname{Im} z)^{a_{1}\gamma_1 p_1'}\mathfrak{g}(\operatorname{Im} u)^{(b_1 - \alpha_1)\gamma_1 p_1' + s_1 p_1' + \alpha_1}}{|Q(z-\overline{u})|^{c_1 \gamma_1 p_1'}} dV(u) \lesssim \mathfrak{g}(\operatorname{Im} z))^{r_1 p_1'},
\]
which, together with (6.3) and (6.4), we further obtain that
\[
\begin{aligned}
& \underset{\eta \in T_{\Lambda_n}}{\mathrm{ess \, sup}} \ \int_{T_{\Lambda_n}} \left[ K_1(z, u) \right]^{\gamma_1 p_1'} \left[ K_2(w, \eta) \right]^{\gamma_2 p_1'} \left[ h_1(u, \eta) \right]^{p_1'} dV_{\alpha_1}(u)\\
& \quad \lesssim \mathfrak{g}(\operatorname{Im} z)^{r_1 p_1'} \mathfrak{g}(\operatorname{Im} w)^{r_2 p_1'} \sim [h_2(z, w)]^{p_1'}.
\end{aligned}
\]

Thus, condition (6.4) holds true for the operator \( S_{\vec{a}, \vec{b}, \vec{c}} \).

We next check the second condition. Observe that
\[
\begin{aligned}
& \int_{T_{\Lambda_n}} [K_1(z, u)]^{\delta_1 q_1} [K_2(w, \eta)]^{\delta_2 q_1} [h_2(z, w)]^{q_1} dV_{\beta_1}(z) \\
& \quad = \frac{ \mathfrak{g}(\operatorname{Im} w)^{a_{2}\delta_2 q_1+\gamma_2 q_1}\mathfrak{g}(\operatorname{Im} u)^{(b_1 - \alpha_1)\mathfrak{g}_1 q_1}}{|Q(w-\overline{\eta})|^{c_2 \mathfrak{g}_2 q_1}} \cdot \int_{T_{\Lambda_n}} \frac{ \mathfrak{g}(\operatorname{Im} z)^{a_{1}\delta_1 q_1+\gamma_1 q_1 + \beta_1}}{|Q(z-\overline{u})|^{c_1 \delta_1 q_1}} dV(z).
\end{aligned}
\]

By the definition of \( r_i \), we obviously have
\begin{equation}
\begin{aligned}
r_i q_i + \beta_i > -1, \end{aligned}
\end{equation}
for any \( i \in \{1, 2\} \). Notice that
\[
(c_i - a_i - b_i + \alpha_i)\delta_i = \tau_i\delta_i = \frac{n + \beta_i}{q_i} + r_i - s_i,
\]
 we have
\begin{equation}
\begin{aligned}
c_i\delta_i q_i - n - 1 - a_i\delta_i q_i - r_i q_i - \beta_i = (b_i - \alpha_i)\delta_i q_i - s_i q_i > 1-\frac{3n}{2}.
\end{aligned}
\end{equation}
Applying this, (6.6), (6.7) and Remark 2.1, we know that,
\[
\int_{T_{\Lambda_n}} \frac{\mathfrak{g}(\operatorname{Im} z)^{\gamma_1 q_1 + \beta_1}}{|Q(z-\overline{u})|^{c_1 \delta_1 q_1}} dv(z) \lesssim \mathfrak{g}(\operatorname{Im} u)^{s_1 q_1 - (b_1 - \alpha_1)\delta_1 q_1}.
\]

From this and Remark 2.2, we deduce that
\[
\int_{T_{\Lambda_n}} \frac{\mathfrak{g}(\operatorname{Im}z))^{a_{1}\delta_1 q_1+\gamma_1 q_1 + \beta_1}}{|Q(z-\overline{u})|^{c_1 \delta_1 q_1}} dV(z) \lesssim \mathfrak{g}(\operatorname{Im} u)^{s_1 q_1 - (a_{1}+b_1 - \alpha_1)\delta_1 q_1}.
\]

This, together with (6.6), (6.7), and Remark 2.2 again, further implies that
\[
\begin{aligned}
& \int_{T_{\Lambda_n}} \left[ \int_{T_{\Lambda_n}} [K_1(z, u)]^{\delta_1 q_1} [K_2(w, \eta)]^{\delta_2 q_1} [h_2(z, w)]^{q_1} dV_{\beta_1}(z) \right]^{q_2 / q_1} dV_{\beta_2}(w) \\
& \quad \lesssim \mathfrak{g}(\operatorname{Im} u)^{s_1 q_2} \int_{T_{\Lambda_n}} \frac{\mathfrak{g}(\operatorname{Im} w))^{\gamma_2 q_2 + \beta_2}}{|Q(w-\overline{\eta})|^{c_2 \delta_2 q_2}} dV(\eta) \\
& \quad \lesssim \mathfrak{g}(\operatorname{Im} u)^{s_1 q_2} \sim [h_1(u, \eta)]^{q_2}.
\end{aligned}
\]

 Thus, the operator \( S_{\vec{a}, \vec{b}, \vec{c}} \) satisfies all conditions of Lemma 2.4.

 \hfill $\square$
\section{Proof of Theorem 5}
We consider the test function \( f_R(z,w) \) defined by
\begin{equation}
f_R(z,w) = \frac{1}{Q(z + \text{iR})^{s_{1}}Q(w + \text{iR})^{s_{2}}} \quad z,w \in T_{\Lambda_n},
\end{equation}

where \( r > 0 \), \( R = (0', r) \in \mathbb{R}^{n - 1} \times \mathbb{R} \), and the real parameters \( s_{i}\) satisfy the conditions
\begin{equation}
\begin{cases}
s_{i} >  \max\left\{ \frac{n}{2}-1, \frac{n-1}{p_{i}} \right\}\\
s_{i} - l_{i} > \max\left\{ \frac{\alpha - 1}{p_{i}}+\frac{3n}{2p}, \frac{3n}{2} - 1 - c_{i} + b_{i} \right\}
\end{cases}.
\end{equation}
We first compute the norm of \( f_R(z,w) \) in $L_{\vec{\alpha}}^{\vec{p}}$. Using Remark 2.1, we deduce that
\begin{equation}
\begin{aligned}
||f_R(z,w)||_{\vec{\alpha}}^{\vec{p}}
&=
\int_{T_{\Lambda_n}}
\int_{T_{\Lambda_n}}
\lvert f(z, w) \rvert \, dV_{\alpha_1}(z)
\, dV_{\alpha_2}(w)
\\
&=
\int_{T_{\Lambda_n}}
\int_{T_{\Lambda_n}}
\lvert \frac{1}{Q(z + \text{iR})^{s_{1}}Q(w + \text{iR})^{s_{2}}} \rvert \, dv_{\alpha_1}(z)
\, dv_{\alpha_2}(w)
\\
&= \int_{T_{\Lambda_n}} \frac{\mathfrak{g}(\operatorname{Im} z)^{\alpha_{1}}}{|Q(z + {iR})|^{s_{1}} } dV(z) \cdot \int_{T_{\Lambda_n}} \frac{\mathfrak{g}(\operatorname{Im} w)^{\alpha_{2}}}{|Q(w + {iR})|^{s_{2}} } dV(w) \\
&=CR^{n-s_{i}+\alpha_{i}}.
\end{aligned}
\end{equation}

where i = 1, 2 and \( C \) is a constant that depends on the parameters \( \alpha_{i}, p_{i},  l_{i}, s_{i} \). The condition (7.2) guarantees that the function \( f_R(z,w) \) belongs to $L_{\vec{\alpha}}^{\vec{p}}$.

We next calculate the norm of \( Tf_R(z,w) \). By Lemma 2.1 and the condition (7.2), we have
$$\begin{aligned}
T_{\vec{a},\vec{b},\vec{c}} f(z,w)
&=
\int_{T_{\Lambda_n}} \int_{T_{\Lambda_n}}
\frac{\mathfrak{g}(\operatorname{Im} z)^{a_1} \mathfrak{g}(\operatorname{Im} \omega)^{a_2}\mathfrak{g}(\operatorname{Im} u)^{b_1} \mathfrak{g}(\operatorname{Im} \eta)^{b_2}}
{Q(z-\overline{u})^{c_1}Q(u+iR)^{s_{1}} Q(\omega-\overline{\eta})^{c_2}Q(\eta+iR)^{s_{2}}}
 dV(u) \, dV(\eta)\\
&=C_{5}\frac{\mathfrak{g}(\operatorname{Im} z)^{a_1} \mathfrak{g}(\operatorname{Im} \omega)^{a_2}}{Q(z+iR)^{c_{1}+s_{1}-n-b_{1}}Q(w+iR)^{c_{2}+s_{2}-n-b_{2}}}.
\end{aligned}
$$

Since the operator \( T \) is bounded on $L_{\vec{\beta}}^{\vec{q}}$, the function \( Tf_R(z,w) \) is in $L_{\vec{\beta}}^{\vec{q}}$. Again by Lemma 2.1 and the condition (7.2), we obtain

\begin{equation}
\begin{cases}
q_{i}a_i + \beta_i > -1 \\
q_{i}(c_{i} - b_{i} - a_{i} - n + s_i) - \beta_i >\frac{3n}{2}-1,
\end{cases} \, i = 1, 2.
\end{equation}

Moreover,
$$
\| T f_{\mathbf{R}} \|_{L_{\vec{\beta}}^{\vec{q}}} = C' R^{a_{i} + b_{i} - c_{i} - s_i + n  + \frac{\beta_1  + n}{q_{1}}},
$$
\( C' \) is a constant depending only on $ a_{i}, b_{i}, c_{i}, l_{1}, s_{i}$, $\vec{\alpha}, \vec{\beta}$, $p_{2}$ and $q_{i}$, where $i$ = 1, 2. Due to the boundedness of the operator \( T \) from $L_{\vec{\alpha}}^{\vec{p}}$ to $L_{\vec{\beta}}^{\vec{q}}$, we have
$$
 \| T f_{\mathbf{R}} \|_{L_{\vec{\beta}}^{\vec{q}}} \leq C''\| f_{\mathbf{R}} \|_{L_{\vec{\alpha}}^{\vec{p}}}.
$$
That is,
$$
C'R^{a_{i} + b_{i} - c_{i} - s_i + n  + \frac{\beta_i  + n}{q_{i}}} \leq CC'' R^{n-s_{i}+\alpha_{i}},
$$
where \( C, C' \) and \( C'' \) are independent of $R$. For the selection of $R$, it is only required to be positive integers. Therefore, for the above formula to hold, the following condition must be satisfied:
\begin{equation}
\begin{cases}
\alpha_{i}<b_{i}\\
c_i = a_i +b_{i}-\alpha_{i}+ \dfrac{n + \beta_i}{q_i}.
\end{cases}
\end{equation}

Combining conditions (7.2) and (7.5), condition (7.4) is equivalent to
\begin{equation}
-a_i q_{i} < \beta_i + 1, \quad i= 1, 2.
\end{equation}

(ii) Proof.  Suppose \( c_i \neq 0 \) for \( i \in \{1,2\} \). From the fact \( -(1 + \beta_i)/q_i < 0 \) for \( i \in \{1,2\} \), we deduce that there exist two negative numbers \( r_1 \) and \( r_2 \) such that, for any \( i \in \{1,2\} \), \( -\frac{1 + \beta_i}{q_i} < r_i < 0 \). For any \( i \in \{1,2\} \), let $\gamma_i :=\frac{s_i-r_{i}}{\tau_i}$  and  \( \delta_i:= \frac{(n + \beta_i)/q_i + r_i - s_i}{\tau_i} \). Obviously, \( \gamma_i + \delta_i = 1 \).
Next, we use Lemma 2.3 to proof.
$$
\begin{aligned}
&[K_1(z,u)]^{\gamma_1} [K_2(w,\eta)]^{\gamma_2} h_1(u,\eta)\\
&= \frac{\mathfrak{g}(\operatorname{Im} z)^{a_{1}\gamma_1}\mathfrak{g}(\operatorname{Im} w)^{a_{2}\gamma_2}(\operatorname{Im} u)^{(b_1 - \alpha_1)\gamma_1 + s_1}\mathfrak{g}(\operatorname{Im} \eta)^{(b_2 - \alpha_2)\gamma_2 + s_2}}{|Q(z-\overline{u})|^{c_1 \gamma_1} |Q(\omega-\overline{\eta})|^{c_2 \gamma_2}}.
\end{aligned}$$

For any \( i \in \{1,2\} \), we get
\begin{equation}
c_i \gamma_i = (b_i - \alpha_i + a_i)\gamma_i + s_i - r_i
\end{equation}
due to \( \tau_i = c_i - a_i - b_i + \alpha_i \) and \( \gamma_i = (s_i - r_i)/\tau_i \).

Then, according to (7.7), for any given \( z \in T_{\Lambda_n}\) and \( u \in T_{\Lambda_n} \), we have
\[
\frac{\mathfrak{g}(\operatorname{Im} z)^{a_1 \gamma_1} \mathfrak{g}(\operatorname{Im} u)^{(b_1 - \alpha_1)\gamma_1 + s_1}}{|Q(z-\overline{u})|^{c_1 \gamma_1}} = \left( \frac{\mathfrak{g}(\operatorname{Im} u)}{|Q(z-\overline{u})|} \right)^{(b_1 - \alpha_1)\gamma_1 + s_1} \left( \frac{\mathfrak{g}(\operatorname{Im} z)}{|Q(z-\overline{u})|} \right)^{a_1 \gamma_1} \lesssim \mathfrak{g}(\operatorname{Im} z)^{r_1}
\]
and, similarly, for any given \( w \in T_{\Lambda_n} \) and any \( \eta \in T_{\Lambda_n} \),
\[
\frac{\mathfrak{g}(\operatorname{Im} w)^{a_2 \gamma_2}\mathfrak{g}(\operatorname{Im} \eta)^{(b_2 - \alpha_2)\gamma_2 + s_2}}{|Q(\omega-\overline{\eta})|^{c_2 \gamma_2}} = \left( \frac{\mathfrak{g}(\operatorname{Im} \eta)}{|Q(\omega-\overline{\eta})|} \right)^{(b_2 - \alpha_2)\gamma_2 + s_2} \left( \frac{\mathfrak{g}(\operatorname{Im} w)}{|Q(\omega-\overline{\eta})|} \right)^{a_2 \gamma_2} \lesssim\mathfrak{g}(\operatorname{Im} w)^{r_2}
\]

Thus, for any given $(z,w)\in T_{\Lambda_n}$
\[
\underset{(u,\eta) \in T_{\Lambda_n} \times T_{\Lambda_n}}{\mathrm{ess \, sup}} \, [K_1(z,u)]^{\gamma_1} [K_2(w,\eta)]^{\gamma_2} h_1(u,\eta) \lesssim h_2(z,w).
\]

Thus, condition (3.4) holds true for the operator $S_{\vec{a},\vec{b},\vec{c}}$.

Note that the condition (2.5) in Lemma 2.3 is the same as the condition (2.3) in Lemma 2.2, thus according to the proof of the second part of Theorem 2, we know that the condition (2.5) still holds. Therefore, the operator \( S_{\vec{a},\vec{b},\vec{c}} \) satisfies all the conditions of Lemma 2.3, then operator $S_{\vec{a},\vec{b},\vec{c}}$ is bounded from \( L_{\vec{\alpha}}^{\vec{1}}(T_{\Lambda_n})\) to \( L_{\vec{\beta}}^{\vec{q}}(T_{\Lambda_n} )\).

 \hfill $\square$

\end{document}